\documentclass[11pt,a4paper]{article}

\usepackage{amsmath,amsfonts,amsthm,amssymb,amscd, appendix, verbatim}
\usepackage{xpatch}
\usepackage[utf8]{inputenc}
\usepackage{setspace, color}
\usepackage{array}
\usepackage{cite}
\usepackage[pagewise]{lineno}
\usepackage{bm}

\newcommand{\caixa}{\hglue15.7cm$\square$\vspace{5mm}}
\newtheorem{theorem}{Theorem}[section]
\newtheorem{corollary}[theorem]{Corollary}

\newtheorem{lemma}[theorem]{Lemma}
\newtheorem{proposition}[theorem]{Proposition}
\theoremstyle{definition}

\newtheorem{remark}[theorem]{Remark}

\makeatletter
   \xpatchcmd{\@thm}{\fontseries\mddefault\upshape}{}{}{} 
\makeatother

\usepackage[margin=1in]{geometry}

\title{Temporal decay rates for weak solutions of the  Navier-Stokes Equations with supercritical fractional dissipation}
\author{Wilberclay G. Melo\footnote{Corresponding author: Departamento de Matem\'atica, Universidade Federal de Sergipe, S\~ao Crist\'ov\~ao SE 49100-000,
  Brazil, e-mail: wilberclay@academico.ufs.br. This author is partially supported by CNPq grant 309880/2021-1.}
  }
\date{}

\begin{document}

\maketitle

\begin{abstract}

In this paper, we establish temporal decay for a weak solution $u(x,t)$ (with initial data $u_0$) of the Navier-Stokes equations with supercritical fractional dissipation $\alpha \in (0,\frac{5}{4})$ in $L^2(\mathbb{R}^3)$ and $\dot{H}^s(\mathbb{R}^3)$ ($s\leq0$). More precisely, we prove that $u$ satisfies the following upper bound:
$$	\|u(t)\|_{2}^2\leq C(1+t)^{-\frac{3-2p}{2\alpha}}, \quad\forall t>0.$$
 This estimate leads us to show the next inequality:
$$	\|u(t)\|_{\dot{H}^{-\delta}}^2\leq C(1+t)^{-\frac{3-2\delta-2p}{2\alpha}}, \quad\forall t>0.$$
These results are obtained by applying standard Fourier Analysis and they hold for $\alpha\in(0,\frac{5}{4}),$ $p\in[-1,\frac{3}{2})$, $\delta\in [0, \frac{3-2p}{2})$ and $u_0\in L^2(\mathbb{R}^3)\cap \mathcal{Y}^p(\mathbb{R}^3)$ (and also $u_0\in L^1(\mathbb{R}^3)$ for $p=-1$ and a certain finite set of values of $\alpha$).
\end{abstract}

\vspace{0.2cm}



\section{Introduction}

\hspace{0.5cm} In this work, we prove  temporal decay rates for weak solutions of the following  incompressible
Navier-Stokes equations:
\begin{equation}\label{micropolar}
\left\{
\begin{array}{l}
u_t\,+\, \:\!(-\Delta)^{\alpha} u
\;\!+\,
u \cdot \nabla u
\,+\,
\nabla \;\!\pi
\;=\;
0, \quad x\in \mathbb{R}^3, \quad t > 0,\\
%
%
%
%
%
%
\mbox{div}\:u \;=\;0, \quad x\in \mathbb{R}^3, \quad t > 0,\\
%
%
u(x,0) \,=\, u_0(x), \quad x\in \mathbb{R}^3,
\end{array}
\right.
\end{equation}
where  $u(x,t)\in \mathbb{R}^3$ denotes the incompressible velocity field
and $\pi(x,t)\in \mathbb{R}$ the hydrostatic pressure (see \cite{artigowilberthyagomanasses,novo,nonlinear,MHD3,Ahn,Shan,artigonata,artigothyago,wu,Jiu} and papers included). The initial data for the velocity field $u_{0}$  is assumed to be divergence free, i.e.,
$\mbox{div}\,u_{0}=  0.$ Let us recall that $(-\Delta)^{\alpha}$  represents the fractional Laplacian (see, for example, \cite{fractional,stein} and references therein for more details), by considering the supercritical case $\alpha\in(0,\frac{5}{4})$ (see \cite{nonlinear} for more details of this definition).  It is important to recall that, by the Spectral Theorem,  $(-\Delta)^{\alpha}$ assumes the diagonal form in the Fourier variable, that is, this is a Fourier multiplier operator with symbol $|\xi|^{2\alpha}$. Physically, (\ref{micropolar})  is the equation that describes the motion of a fluid with internal
friction interaction and such motion is a chain of particles that are connected  by elastic springs (we cite \cite{fractional} and references therein for more information).


One of the most relevant open problems in Analysis is related to the standard Navier-Stokes equations, which are   a particular case of the system (\ref{micropolar}). More specifically, when $\alpha=1$, we obtain the famous Navier-Stokes equations below:
\begin{equation}\label{ns}
\left\{
\begin{array}{l}
u_t
\;\!+\,
u \cdot \nabla u
\,+\,
\nabla \;\!\pi
\;=\;
\Delta u, \quad x\in \mathbb{R}^3, \quad t > 0,\\
%
%
%
%
%
%
\mbox{div}\:u \;=\;0, \quad x\in \mathbb{R}^3, \quad t > 0,\\
%
%
u(x,0) \,=\, u_0(x), \quad x\in \mathbb{R}^3,
\end{array}
\right.
\end{equation}
The well-posedness of these equations relies on proving that  singularities for their solutions  are not detected in finite time
by assuming  smooth initial data with finite energy (see \cite{Le34}). Therefore, the  fractional Laplacian in (\ref{micropolar}) may be very helpful to comprehend better how to study (\ref{ns}). More specifically, although the literature does not seem to be ready to show the existence of a global  classical solution for (\ref{micropolar}) whether  $\alpha \in (0,\frac{5}{4})$, J. Wu \cite{wu} proved  that there is a unique global classical solution for the Navier-Stokes equations (\ref{micropolar}) if $\alpha\geq\frac{5}{4}$, $u_0\in H^s(\mathbb{R}^3)$ ($s>2\alpha$).
 On the other hand, it is well known that there exist weak solutions for this same system (\ref{micropolar}) (see \cite{nonlinear,wu} and references therein) such that
\begin{align}\label{*12}
\frac{d}{dt}\|u(t)\|_{2}^2+2\|(-\Delta)^{\frac{\alpha}{2}}u(t)\|_2^2=0
\end{align}
and also
\begin{align}\label{*121}
\|u(t)\|_{2}^2+2\int_0^t\|(-\Delta)^{\frac{\alpha}{2}}u(\tau)\|_2^2d\tau\leq \|u_0\|_{2}^2,\quad\forall t>0,
\end{align}
where $u_0\in L^2(\mathbb{R}^3)$ (see Definition 1.1 and Proposition 3.1 in \cite{nonlinear}). Notice that (\ref{*121}) implies that $u\in L^\infty([0,\infty);$ $L^2(\mathbb{R}^3))$. These statements above are one of the reasons why we are interested in working with  weak solutions of (\ref{micropolar}).

The main motivations for our work are the publications of some papers that study  important results involving temporal decay rates for weak solutions of (\ref{micropolar}) and (\ref{ns})   in $L^2(\mathbb{R}^3)$ and $\dot{H}^{s}(\mathbb{R}^3)$ ($s\leq0$) (see \cite{Niche,Jiu,schonbek,novo,nonlinear,Ahn,Shan,artigowilberthyagomanasses} and references therein). More recently, L. Deng and H. Shang \cite{nonlinear} proved that a weak solution of (\ref{micropolar}) satisfies
\begin{align}\label{nonlinear1}
	\|u(t)\|_{2}^2\leq C(1+t)^{-\frac{3}{2\alpha}},\quad \forall t>0,
\end{align}
where $\alpha\in(0,\frac{5}{4})$ and $u_0\in L^2(\mathbb{R}^3)\cap L^1(\mathbb{R}^3)$. Moreover, \cite{nonlinear} also showed that if the initial data $u_0$ is small enough    in $H^s(\mathbb{R}^3)$ ($s>\frac{5}{2}-\alpha$); then, there is a unique global solution for (\ref{micropolar}) such that, for any $T>0$,
$$u\in L^\infty([0,T];H^s(\mathbb{R}^3))\cap L^2([0,T];H^{s+\alpha}(\mathbb{R}^3))$$
and
\begin{align*}
	\|u(t)\|_{\dot{H}^{-\delta}}\leq C, \quad\forall t>0,
\end{align*}
where $\alpha\in(0,\frac{5}{4})$ and $\delta\in[0,\frac{3}{2})$. Furthermore, by considering the system (\ref{ns}), H. Bae, J. Jung and J. Shin \cite{novo} established the following inequalities:
\begin{align*}
	\|u(t)\|_{2}^2\leq C(1+t)^{-\frac{3-2p}{2}}
\end{align*}
and also
\begin{align*}
	\|u(t)\|_{\dot{H}^{-\delta}}^2\leq C(1+t)^{-\frac{3-2\delta-2p}{2}}, \quad\forall t>0,
\end{align*}
where $p\in [-1,1]$, $\delta\in(0,\frac{3-2p}{2})$ and $u_0\in L^2(\mathbb{R}^3)\cap \mathcal{Y}^p(\mathbb{R}^3)$ (see (\ref{espaconovo}) for definitions and notations).

Now, let us present our main results  related to weak solutions for the Navier-Stokes equations (\ref{micropolar}) and (\ref{ns}). The first one shows decay rates for these solutions in the specific Lebesgue space $L^2(\mathbb{R}^3)$.

\begin{theorem}\label{teoremaexistenciaw1}
Let $\alpha\in(0,\frac{5}{4})$, $p\in [-1,\frac{3}{2})$ and $u_0\in L^2(\mathbb{R}^3)\cap \mathcal{Y}^p(\mathbb{R}^3)$. Assume that $u\in L^\infty([0,\infty);$ $L^2(\mathbb{R}^3))$ is a weak solution for the Navier-Stokes equations \emph{(\ref{micropolar})}. Then, if $p\in[-1,2\alpha-1]$, one has
\begin{align}\label{wn10}
u\in L^\infty([0,\infty);\mathcal{Y}^p(\mathbb{R}^3))
\end{align}
and, if $p\in[-1,\tfrac{3}{2})$, it holds
\begin{align}\label{l2}
	\|u(t)\|_{2}^2\leq C(1+t)^{-\frac{3-2p}{2\alpha}}, \quad\forall  t>0,
\end{align}
where $C$ is a positive constant. Here, we consider also that $u_0\in L^1(\mathbb{R}^3)$ if
\begin{align}\label{condicaoalphan}
p=-1 \,\,\hbox{ and } \,\,\alpha=\frac{5\cdot 2^n-5}{2^{n+2}-2}, \quad  n\in \{1,2,...,m_0+1\},
\end{align}
where $m_0$ is the smallest natural number  such that $m_0>\log_2(\frac{5-2\alpha}{10-8\alpha})-1.$
\end{theorem}

As a direct  consequence of Theorem \ref{teoremaexistenciaw1}, we obtain our second main result.

\begin{corollary}\label{corolario}
  Let $\alpha\in(0,\frac{5}{4})$, $p\in [-1,\frac{3}{2})$, $\delta\in [0, \frac{3-2p}{2})$ and $u_0\in L^2(\mathbb{R}^3)\cap \mathcal{Y}^p(\mathbb{R}^3)$. Assume that $u\in L^\infty([0,\infty);L^2(\mathbb{R}^3))$ is a weak solution for the Navier-Stokes equations \emph{(\ref{micropolar})}. Then,
  \begin{enumerate}
    \item[\emph{i)}] If $p\in[-1,2\alpha-1]$, we have 
    \begin{align}\label{hdelta1}
	\|u(t)\|_{\dot{H}^{-\delta}}^2\leq C(1+t)^{-\frac{3-2\delta-2p}{2\alpha}}, \quad\forall t>0,
\end{align}
where $C$ is a positive constant. Here, we consider also that $u_0\in L^1(\mathbb{R}^3)$ if
\begin{align*}
p=-1 \,\,\hbox{ and } \,\,\alpha=\frac{5\cdot 2^n-5}{2^{n+2}-2}, \quad n\in \{1,2,...,m_0+1\},
\end{align*}
where $m_0$ is the smallest natural number  such that $m_0>\log_2(\frac{5-2\alpha}{10-8\alpha})-1.$
    \item[\emph{ii)}] If $p\in[2\alpha-1,\frac{3}{2})$, we obtain
    \begin{align}\label{hdelta}
	\|u(t)\|_{\dot{H}^{-\delta}}^2&\leq C[(1+t)^{-\frac{3-2\delta-2p}{2\alpha}}+(1+t)^{-\frac{5-4\alpha-2\delta}{2\alpha}}],
\end{align}
where $C$ is a positive constant.
  \end{enumerate}

\end{corollary}

\begin{remark}
It is relevant to emphasize that the upper bounds given by (\ref{hdelta1}) and (\ref{hdelta}) imply that any weak solution $u$ for the Navier-Stokes equations (\ref{micropolar}) must verify the following:
\begin{align}\label{hdelta2}
\|u(t)\|_{\dot{H}^{-\delta}}^2\leq C (1+t)^{-\frac{3-2\delta-2p}{2\alpha}}, \quad \forall t>0,
\end{align}
provided that $\alpha\in(0,\frac{5}{4})$, $p\in[-1,\frac{3}{2})$ and $\delta\in [0, \frac{3-2p}{2})$. 
\end{remark}

When we assume  $\alpha=1$ in Theorem \ref{teoremaexistenciaw1} and Corollary \ref{corolario}, it means that we are studying the Navier-Stokes equations (\ref{ns}). This condition was considered by H. Bae, J. Jung and J. Shin in \cite{novo} and the conclusions obtained by these researchers are given in our main results as well (with more possibilities of choice for $p$, since $p\in[-1,\frac{3}{2})$ instead of $p\in[-1,1]$). Thus, more precisely, our paper generalizes and  improves (since $1\neq (5\cdot 2^n-5)(2^{n+2}-2)^{-1}$ for all $n\in \mathbb{N}$, see (\ref{condicaoalphan})) Theorem 1.1 and Corollary 1.1 established in \cite{novo}.

L. Deng and H. Shang \cite{nonlinear} proved that a weak solution $u$ for the Navier-Stokes equations (\ref{micropolar}) must satisfy the  temporal decay given by (\ref{nonlinear1}),
where  $u_0\in L^2(\mathbb{R}^3)\cap L^1(\mathbb{R}^3)$. It is easy to check that this estimate (\ref{nonlinear1}) is presented in our Theorem \ref{teoremaexistenciaw1} provided that $u_0\in L^2(\mathbb{R}^3)\cap \mathcal{Y}^p(\mathbb{R}^3)$; in fact, it is enough to take $p\in(-1,0]$ in (\ref{l2}). Furthermore, if $p\in(-1,0]$ and $\delta\in[0,\frac{3-2p}{2})$ (notice that  $[0,\frac{3}{2})\subseteq [0,\frac{3-2p}{2})$ whether $p\leq0$), then
\begin{align}\label{hdeltaconstante}
	\|u(t)\|_{\dot{H}^{-\delta}}\leq C, \quad\forall t>0,
\end{align}
is an immediate consequence of the temporal decay (\ref{hdelta2}). This inequality (\ref{hdeltaconstante}) is presented by \cite{nonlinear} (see (1.17) in Theorem 1.6 (ii) (1)). On the other hand, it is worth to point out that by assuming $p=0$, we have that $u_0\in L^2(\mathbb{R}^3)\cap L^1(\mathbb{R}^3)$ implies that
$u_0\in L^2(\mathbb{R}^3)\cap \mathcal{Y}^0(\mathbb{R}^3)$. Thereby, we also present a greater variety of  possibilities for the initial data $u_0$ in our main results than  Theorem 1.2 ii) obtained by \cite{nonlinear}, for instance.


\section{Prelude}\label{preliminares}

\hspace{0.5cm} This section  presents some definitions, notations and preliminary results that will play an important role in our paper.\\\\
$\bullet$ \textbf{Basic notations and definitions}:\\\\
Let us list the main notations and definitions of this work.
\begin{enumerate}
\item $S'(\mathbb{R}^3)$ is the space of tempered distributions.

\item The Fourier transform and its inverse are defined by
 $$\mathcal{F}(f)(\xi)= \widehat{f}(\xi):=\int_{\mathbb{R}^3}e^{-i\xi\cdot x}f(x)\,dx,$$
and $$\mathcal{F}^{-1}(g)(x):=(2\pi)^{-3}\int_{\mathbb{R}^3}e^{i\xi\cdot x}g(\xi)\,d \xi,$$
repectively.


\item The fractional Laplacian $(-\Delta)^{\alpha}$ (for more details, see \cite{stein}), $\alpha\in(0,\frac{5}{4})$, is defined by
$$\mathcal{F}[(-\Delta)^{\alpha}f](\xi)=|\xi|^{2\alpha} \hat{f} (\xi),\quad\forall \xi\in \mathbb{R}^3,$$
where $f\in S'(\mathbb{R}^3)$ and $\hat{f}\in L^1_{loc}(\mathbb{R}^3).$
\item  The tensor product is given by
$$f\otimes g:=(g_{1}f,g_{2}f,g_{3}f),$$
where $f=(f_1,f_2,f_3)$ and $g=(g_1,g_2,g_3)\in S'(\mathbb{R}^3).$

\item Let $p\in [-1,\frac{3}{2})$. We define (see \cite{Bh} for more details)
\begin{align}\label{espaconovo}
		\mathcal{Y}^p(\mathbb{R}^3):=\{f\in S'(\mathbb{R}^3): \hat{f}\in L^1_{loc}(\mathbb{R}^3) \hbox{  and  } \sup_{\xi\in \mathbb{R}^3}\{|\xi|^p|\widehat{f}(\xi)|\}<\infty\}
\end{align}
 and $\mathcal{Y}^p(\mathbb{R}^3)$-norm is given by
$$ \|f\|_{\mathcal{Y}^p}=\sup_{\xi\in \mathbb{R}^3}\{|\xi|^p|\widehat{f}(\xi)|\}.$$

  \item Let   ($X,\|\cdot\|_{X}$) be a normed space. We define
  $$L^\infty([0,\infty);X)=\{f:[0,\infty)\rightarrow X \hbox{ mensurable function}: \sup_{t\in [0,\infty)}\{\|f(t)\|_{X}\}<\infty\},$$ and $L^{\infty}([0,\infty);X)$-norm is given by
$$\|f\|_{L^{\infty}([0,\infty);X)}:=\sup_{t\in [0,\infty)}\{\|f(t)\|_{X}\}.$$

\item The constants in this paper may change their values from line to line without change of notation.
\end{enumerate}



\hspace{-0.6cm}$\bullet$ \textbf{Auxiliary results}:\\\\
The next result extends Lemma 2.1 established in \cite{novo} in order  to be applicable to the case of the Navier-Stokes equations (\ref{micropolar}).

\begin{lemma}\label{lemaderivada}
Assume that $\alpha>0$ and $f$ is a smooth function such that $f(t)\in L^2(\mathbb{R}^3)$, for all $t\geq0$, and
\begin{align}\label{lemaderivada1}
\frac{d}{dt}\|f(t)\|_{2}^2+ \|(-\Delta)^\frac{\alpha}{2} f(t)\|_{2}^2\leq 0,\quad\forall t\geq 0.
\end{align}
Suppose that there are  constants  $C_1>0$ and $p<\frac{3}{2}$ that satisfy the following inequality:
\begin{align}\label{hipoteselema2.1}
\sup_{t\in[0,\infty)}\{\sup_{|\xi|\leq 1} \{|\xi|^p|\hat{f}(\xi,t)|\}\}\leq C_1.
\end{align}
Then, we obtain
\begin{align*}
\|f(t)\|_{2}^2\leq C(1+t)^{-\frac{3-2p}{2\alpha}},\quad \forall t>0,
\end{align*}
where $C$ is a positive constant.
\end{lemma}
\begin{proof}
First of all, let $q$ be a real number such that $q>\max\{1,\frac{3-2p}{2\alpha}\}$. Secondly, by applying (\ref{lemaderivada1}) and Plancherel's identity, it is true that
\begin{align*}
\frac{d}{dt}\|f(t)\|_{2}^2&\leq  -C\int_{\{|\xi|> (\frac{q}{1+t})^\frac{1}{2\alpha}\}}|\xi|^{2\alpha}| \hat{f}(\xi)|^2d\xi \\
&\leq-\frac{Cq}{1+t}\int_{\{|\xi|> (\frac{q}{1+t})^\frac{1}{2\alpha}\}} |\hat{f}(\xi)|^2d\xi\\
&= -\frac{q}{1+t}\|f(t)\|_{2}^2+\frac{Cq}{1+t}\int_{\{|\xi|\leq (\frac{q}{1+t})^\frac{1}{2\alpha}\}}| \hat{f}(\xi)|^2d\xi,
\end{align*}
for all $t\geq0$ (since $\alpha>0$). Consequently, (\ref{hipoteselema2.1}) implies that
\begin{align*}
\frac{d}{dt}\|f(t)\|_{2}^2 +\frac{q}{1+t}\|f(t)\|_{2}^2&\leq \frac{C}{1+t}\int_{\{|\xi|\leq (\frac{q}{1+t})^\frac{1}{2\alpha}\}} |\xi|^{-2p}d\xi\\
&=\frac{C}{1+t}\int_0^{(\frac{q}{1+t})^\frac{1}{2\alpha}} r^{2-2p}dr\\
&= C (1+t)^{-1-\frac{3-2p}{2\alpha}},
\end{align*}
for all $t\geq q-1>0$ (recall that $q>1$, $3-2p>0$ and $\alpha>0$).  Multiplying the inequality above by $(1+t)^q$, one has
\begin{align*}
(1+t)^q\frac{d}{dt}\|f(t)\|_{2}^2 +q(1+t)^{q-1}\|f(t)\|_{2}^2&\leq C (1+t)^{q-1-\frac{3-2p}{2\alpha}},\quad \forall t\geq q-1.
\end{align*}
 Thereby, we conclude that
\begin{align*}
\frac{d}{dt}\{(1+t)^q\|f(t)\|_{2}^2\}&\leq C (1+t)^{q-1-\frac{3-2p}{2\alpha}},\quad\forall t\geq q-1.
\end{align*}
Now, integrate this last inequality over the interval $[q-1,t]$ (with $t\geq q-1>0$) in order to obtain
\begin{align*}
(1+t)^q\|f(t)\|_{2}^2&\leq q^q\|f(q-1)\|_{2}^2+ C \int_{q-1}^t(1+\tau)^{q-1-\frac{3-2p}{2\alpha}}d\tau,\quad \forall t\geq q-1.
\end{align*}
As a result, we deduce that
\begin{align*}
\|f(t)\|_{2}^2&\leq C[(1+t)^{-q}+  (1+t)^{-\frac{3-2p}{2\alpha}}],\quad \forall t\geq q-1,
\end{align*}
since $q-\frac{3-2p}{2\alpha}>0.$ Therefore, it follows that
\begin{align}\label{wil1}
\|f(t)\|_{2}^2&\leq C(1+t)^{-\frac{3-2p}{2\alpha}},\quad \forall t\geq q-1,
\end{align}
once again because $q>\frac{3-2p}{2\alpha}.$

On the other hand, (\ref{lemaderivada1}) implies that
\begin{align*}
\|f(t)\|_{2}^2+ \int_0^t\|(-\Delta)^\frac{\alpha}{2} f(\tau)\|_{2}^2d\tau\leq \|f(0)\|_{2}^2,\quad\forall t> 0.
\end{align*}
Then, we have
\begin{align*}
\|f(t)\|_{2}^2\leq \|f(0)\|_{2}^2 ,\quad\forall t>0.
\end{align*}
This leads us to conclude that
\begin{align*}
(1+t)^{\frac{3-2p}{2\alpha}}\|f(t)\|_{2}^2\leq q^{\frac{3-2p}{2\alpha}}\|f(0)\|_{2}^2,\quad\forall t\in(0,q-1),
\end{align*}
since $3-2p>0$ and $\alpha>0.$ As a consequence, it holds
\begin{align}\label{wil2}
\|f(t)\|_{2}^2\leq C(1+t)^{-\frac{3-2p}{2\alpha}},\quad\forall t\in(0,q-1).
\end{align}

From (\ref{wil1}) and (\ref{wil2}), one obtains
\begin{align*}
\|f(t)\|_{2}^2\leq C(1+t)^{-\frac{3-2p}{2\alpha}},\quad\forall t>0.
\end{align*}
This proves Lemma \ref{lemaderivada}.

\end{proof}

Our preliminary result below establishes necessary conditions in order to obtain the  decay rates presented in Theorem \ref{teoremaexistenciaw1}  for a weak solution $u$ of the Navier-Stokes equations (\ref{micropolar}) in $L^2(\mathbb{R}^3)$. Furthermore, the next lemma  is a useful tool to show (\ref{wn10}) and that
        \begin{align*}
    \sup_{t\in[0,\infty)}\{\sup_{|\xi|\leq 1} \{|\xi|^p|\hat{u}(\xi,t)|\}\}\leq C,\quad\forall p\geq2\alpha-1.
    \end{align*}

\begin{lemma}\label{lema2.1}
Let $\alpha\in(0,\frac{5}{4})$ and define the following sequences:
\begin{enumerate}
\item[\emph{1)}] $s_n=2^{n+1}-2$, for all $n\in \mathbb{N}\cup \{0\};$
\item[\emph{2)}] $p_0=2\alpha-1$ and $p_n=6\alpha-6+s_{n-1}(4\alpha-5)$, for all $n\in \mathbb{N};$
\item[\emph{3)}] $\alpha_0=0$ and $\displaystyle\alpha_n=\frac{5(1+s_{n-1})}{2(3+2s_{n-1})}$, for all $n\in \mathbb{N}.$
\end{enumerate}
Then,  $(s_n)_{n\geq0}$ and $(\alpha_n)_{n\geq0}$ are increasing sequences and $(p_n)_{n\geq0}$ is a decreasing sequence. Moreover, it is true that
\begin{enumerate}
\item[\emph{i)}] $ s_n\geq0$, for all $n\in \mathbb{N}\cup\{0\};$
\item[\emph{ii)}] $ s_n=2+2s_{n-1}$, for all $n\in \mathbb{N};$
\item[\emph{iii)}] $p_{n}=2\alpha-4+2p_{n-1}$, for all $n\in \mathbb{N};$
\item[\emph{iv)}] $p_n<\frac{3}{2}$, for all $n\in \mathbb{N}\cup\{0\};$
\item[\emph{v)}] $0<\alpha_n<\frac{5}{4}$, for all $n\in \mathbb{N};$
\item[\emph{vi)}] $\alpha>\alpha_n\Leftrightarrow p_n>-1$, for all $n\in \mathbb{N}\cup\{0\};$
\item[\emph{vii)}] $\alpha<\alpha_n\Leftrightarrow p_n<-1$, for all $n\in \mathbb{N};$
\item[\emph{viii)}] $\alpha=\alpha_n\Leftrightarrow p_n=-1$, for all $n\in \mathbb{N}.$
\end{enumerate}
\end{lemma}

\begin{proof}
First of all, by 1), note that
\begin{align}\label{wf1}
s_{n+1}=2^{n+2}-2> 2^{n+1}-2=s_n,\quad\forall n\in \mathbb{N}\cup \{0\}.
\end{align}
This shows that $(s_n)_{n\geq0} $ is a increasing sequence.

Consequently, from 2), one concludes
\begin{align}\label{wf3}
p_{n+1}<p_n \Leftrightarrow s_{n}(4\alpha-5)<s_{n-1}(4\alpha-5)\Leftrightarrow s_n> s_{n-1},
\end{align}
for all $n\in \mathbb{N}$ (since $4\alpha-5<0$). In addition, by 2), we have
\begin{align}\label{wf2}
p_{1}<p_0 \Leftrightarrow 6\alpha-6<2\alpha-1\Leftrightarrow \alpha<\frac{5}{4}.
\end{align}
Thereby, (\ref{wf1}), (\ref{wf3}) and (\ref{wf2}) infer that $(p_n)_{n\geq0}$ is a decreasing sequence.

It is easy to check that 1) implies i), since that
\begin{align*}
2^{n+1}\geq 2,\quad\forall n\in \mathbb{N}\cup\{0\}.
\end{align*}

Notice that 3) and i) lead us to obtain
\begin{align}\label{wf4}
\alpha_{n}<\alpha_{n+1} \Leftrightarrow \frac{1+s_{n-1}}{3+2s_{n-1}}<\frac{1+s_{n}}{3+2s_{n}}\Leftrightarrow s_{n-1}< s_{n},
\end{align}
for all $n\in \mathbb{N}$. In addition, by 1) and 3), we have
\begin{align}\label{wf5}
\alpha_{0}<\alpha_1 \Leftrightarrow 0<\frac{5}{6}.
\end{align}
Thereby, from (\ref{wf1}), (\ref{wf4}) and (\ref{wf5}), it results that $(\alpha_n)_{n\geq0}$ is a increasing sequence.

It is easy to see that the definition given in 1) implies the following equalities:
\begin{align}\label{wil3}
2+2s_{n-1}=2^{n+1}-2=s_n,\quad \forall n\in \mathbb{N}.
\end{align}
This proves ii).

Notice also that, by (\ref{wil3}) and 2), we have
\begin{align*}
2\alpha-4+2p_{n-1}&=6\alpha-6+(2+2s_{n-2})(4\alpha-5)\\
&= 6\alpha-6+s_{n-1}(4\alpha-5)=p_{n},
\end{align*}
for all $n\geq2$. Furthermore, it follows that
\begin{align*}
2\alpha-4+2p_{0}&=6\alpha-6=p_{1},
\end{align*}
see 1) and 2). Thus, iii) has been proved.

In order to verify iv), we observe that
\begin{align*}
p_n= 6\alpha-6+s_{n-1}(4\alpha-5)\leq 6\alpha-6<\frac{3}{2},
\end{align*}
because of the fact that $\alpha\in (0,\frac{5}{4})$ and $s_{n-1}\geq0$ (by i)), for all $n\in \mathbb{N}$. Moreover,
$p_0=2\alpha-1<\frac{3}{2},$ by using again the condition   $\alpha\in (0,\frac{5}{4})$.

As an immediate consequence of i) and 3), one infers that  $\alpha_n>0$, for all $n\in \mathbb{N}$. In addition, we can write, by 3) and i) once more, that
\begin{align*}
\alpha_n<\frac{5}{4}\Leftrightarrow \frac{5(1+s_{n-1})}{2(3+2s_{n-1})}<\frac{5}{4}\Leftrightarrow 2+2s_{n-1}<3+2s_{n-1}\Leftrightarrow 2<3,
\end{align*}
for all $n\in \mathbb{N}$. These last arguments establish v).

To finish this proof, we shall show vi). More precisely, one has
\begin{align*}
p_n>-1\Leftrightarrow 6\alpha-6+s_{n-1}(4\alpha-5)>-1\Leftrightarrow \alpha>\frac{5+5s_{n-1}}{6+4s_{n-1}}=\alpha_n,
\end{align*}
for all $n\in \mathbb{N}$ (see 2), 3) and i)). It is also true that
\begin{align*}
p_0>-1\Leftrightarrow 2\alpha-1>-1\Leftrightarrow \alpha>0=\alpha_0,
\end{align*}
by applying 2) and 3).

The items vii) and viii) are analogous to vi).
\end{proof}

\begin{remark}
By passing to limit in the definition 1) of Lemma \ref{lema2.1}, as $n\rightarrow\infty$, one reaches
\begin{align*}
\lim_{n\rightarrow\infty} s_n = \lim_{n\rightarrow\infty} (2^{n+1}-2)=\infty.
\end{align*}
Consequently, there is $m_0\in \mathbb{N}$ (take the smallest one) such that
\begin{align}\label{sm0}
s_{m_0}> \frac{14\alpha-15}{2(5-4\alpha)}.
\end{align}
This inequality (\ref{sm0}) is equivalent to
\begin{align*}
\alpha< \frac{15+10s_{m_0}}{14+8s_{m_0}}=\frac{5+5(2+2s_{m_0})}{6+4(2+2s_{m_0})}= \frac{5+5s_{m_0+1}}{6+4s_{m_0+1}}=\alpha_{m_0+2},
\end{align*}
by Lemma \ref{lema2.1} i), ii) and 3) (recall that $\alpha\in(0,\frac{5}{4})$), that is,
\begin{align}\label{alphaepm0}
\alpha< \alpha_{m_0+2}.
\end{align}
\end{remark}

Now, let us finish this section by recalling two elementary results.

\begin{lemma}\label{lemaexponencial}
The following statements hold:
\begin{enumerate}
\item[\emph{i)}] Let $a, b > 0$. Then, we have
$$\lambda^ae^{-b\lambda}\leq a^a(eb)^{-a},\quad\forall \lambda>0;$$
\item[\emph{ii)}] Let $a, b > 0$ such that $a+b=1$. Thereby, we obtain
$$\int_0^t(t-\tau)^{-a}\tau^{-b}d\tau=\int_0^1(1-y)^{-a}y^{-b}dy=:\beta(a;b),$$
where $\beta$ is the standard beta function.
\end{enumerate}
\end{lemma}
\begin{proof}
It is enough to apply  Calculus.
\end{proof}

\section{Proof of our main results:}\label{secaoteoremaexistenciaB}

\hspace{0.5cm} Let us  present the proofs of Theorem \ref{teoremaexistenciaw1} and Corollary \ref{corolario}. Thus,  in order to establish the veracity of this first result, it is necessary to show a proposition that plays an important role in this work.

\begin{proposition}\label{inducao}
  Assume that $\alpha\in(0,\frac{5}{4})$, $p\in \mathbb{R}$   and $u_0\in L^2(\mathbb{R}^3)\cap \mathcal{Y}^p(\mathbb{R}^3)$. Let $u\in L^\infty([0,\infty);$ $L^2(\mathbb{R}^3))$ be a weak solution of the Navier-Stokes equations \emph{(\ref{micropolar})}.  Then, for each $n\in \mathbb{N}\cup \{0\}$, we infer
  \begin{enumerate}
    \item[\emph{i)}] If $p\leq p_n$ and $\alpha\in(\alpha_n,\frac{5}{4})$; hence,
    \begin{align*}
    \|u(t)\|_2^2\leq C(1+t)^{-\frac{3-2p_n}{2\alpha}},\quad \forall t>0;
    \end{align*}
        \item[\emph{ii)}] If $p_{n+1}\leq p< p_n$, $p>-1$ and $\alpha\in(\alpha_n,\frac{5}{4})$; then,
        \begin{align*}
    \|u\|_{L^\infty([0,\infty);\mathcal{Y}^p)}\leq \|u_0\|_{\mathcal{Y}^p} + C\,\beta(\tfrac{p+1}{2\alpha};\tfrac{2\alpha-1-p}{2\alpha});
    \end{align*}
    \item[\emph{iii)}] If $\alpha_{n}< \alpha< \alpha_{n+1}$ and $p=-1$; thus,
        \begin{align*}
    \|u\|_{L^\infty([0,\infty);\mathcal{Y}^{-1})}\leq \|u_0\|_{\mathcal{Y}^{-1}} +  \frac{2\alpha C}{-(p_{n+1}+1)}
    \end{align*}
    \end{enumerate}
   where $C$ stands for a positive constant and $\beta$ is the standard beta function \emph{(}see Lemma \emph{\ref{lemaexponencial} ii))}.

   Furthermore, if $p\geq2\alpha-1$, the inequality below holds:
        \begin{align}\label{2alpha-1}
    \sup_{t\in[0,\infty)}\{\sup_{|\xi|\leq 1} \{|\xi|^p|\hat{u}(\xi,t)|\}\}\leq \|u_0\|_{\mathcal{Y}^{p}} + \|u_0\|_2^2
    \end{align}
    and, in particular, if $p=2\alpha-1$, we deduce
    \begin{align}\label{wn}
    \|u\|_{L^\infty([0,\infty);\mathcal{Y}^{2\alpha-1})}\leq \|u_0\|_{\mathcal{Y}^{2\alpha-1}} + \|u_0\|_2^2,
    \end{align}
   \emph{(}See Lemma \emph{\ref{lema2.1}} \emph{2)} and \emph{3)} for the precise definitions of $\alpha_n$ and $p_n$\emph{)}.
\end{proposition}

\begin{proof}
First of all, it is necessary to  apply the heat semigroup $e^{-  (t-\tau)(-\Delta)^{\alpha}}$ (with $\tau\in[0,t]$) to the first equation in (\ref{micropolar}) to deduce
\begin{align}\label{livro1}
e^{-  (t-\tau)(-\Delta)^{\alpha}}u_\tau+
e^{-  (t-\tau)(-\Delta)^{\alpha}}P(u \cdot \nabla u)
+
e^{-  (t-\tau)(-\Delta)^{\alpha}}(-\Delta)^{\alpha} u
=
0.
\end{align}
where $P$ is Leray's projector. It is known that this operator satisfies the following:
\begin{align}\label{transformadadeP}
|\mathcal{F}[P(f)](\xi)|\leq|\widehat{f}(\xi)|,\quad\forall \xi\in \mathbb{R}^3.
\end{align}
Integrate (\ref{livro1}) over $[0,t]$ to establish the equation below:
\begin{align}\label{novo3}
u(t)= e^{- t(-\Delta)^{\alpha}}u_0 - \int_{0}^te^{-  (t-\tau)(-\Delta)^{\alpha}} P(u\cdot\nabla u)(\tau)d\tau .
\end{align}
By using the Fourier transform in the equality (\ref{novo3}), applying (\ref{transformadadeP}), Plancherel's identity and Young's inequality, we obtain
\begin{align}\label{wil5}
\nonumber|\widehat{u}(\xi,t)|&\leq e^{- t|\xi|^{2\alpha}}|\widehat{u}_0(\xi)| + \int_{0}^t |\xi|e^{-  (t-\tau)|\xi|^{2\alpha}} |\mathcal{F}[u\otimes u](\tau)|d\tau\\
&\leq |\widehat{u}_0(\xi)| + \int_{0}^t |\xi|e^{-  (t-\tau)|\xi|^{2\alpha}} \|u(\tau)\|_2^2d\tau.
\end{align}
We are now ready to prove Proposition \ref{inducao} i)--iii) through an inductive process related to $n\in \mathbb{N}\cup\{0\}$. Thus, let us start with the case $n=0$.\\\\
\underline{1º Case}: Assume $n=0.$\\\\
i) Consider that $p\leq 2\alpha-1$ and $\alpha\in(0,\frac{5}{4})$ (recall that $\alpha_0=0$ and $p_0=2\alpha-1$, see Lemma \ref{lema2.1} 2) and 3)). Thus,  multiply the inequality (\ref{wil5}) by $|\xi|^{2\alpha-1}$ and apply (\ref{*121}) to obtain
\begin{align*}
\nonumber|\xi|^{2\alpha-1}|\widehat{u}(\xi,t)|&\leq |\xi|^{2\alpha-1-p}|\xi|^{p}|\widehat{u}_0(\xi)| + \int_{0}^t |\xi|^{2\alpha}e^{-  (t-\tau)|\xi|^{2\alpha}} \|u(\tau)\|_2^2d\tau\\
&\leq |\xi|^{p}|\widehat{u}_0(\xi)| + \|u_0\|_2^2\int_{0}^t |\xi|^{2\alpha}e^{-  (t-\tau)|\xi|^{2\alpha}} d\tau\\
&\leq \|u_0\|_{\mathcal{Y}^p} + \|u_0\|_2^2,
\end{align*}
for all $\xi\in \mathbb{R}^3$ such that $|\xi|\leq1$ and $t\geq0$ (since $2\alpha-1-p\geq0$). As a result, we deduce that
\begin{align*}
\sup_{t\in [0,\infty)}\{\sup_{|\xi|\leq 1} |\xi|^{2\alpha-1}|\widehat{u}(\xi,t)|\}\leq \|u_0\|_{\mathcal{Y}^p} + \|u_0\|_2^2.
\end{align*}
The assumption  $u_0\in L^2(\mathbb{R}^3)\cap \mathcal{Y}^p(\mathbb{R}^3)$, the fact that $2\alpha-1<\frac{3}{2}$ (because $\alpha\in(0,\frac{5}{4})$), (\ref{*12}) and Lemma \ref{lemaderivada} imply that
\begin{align}\label{P1}
    \|u(t)\|_2^2\leq C(1+t)^{-\frac{3-2p_0}{2\alpha}},\quad \forall t>0,
    \end{align}
    where $p_0=2\alpha-1$ (see Lemma \ref{lema2.1} 2)). This proves i) with $n=0.$\\\\
ii) Consider that $6\alpha-6\leq p< 2\alpha-1$, $p>-1$ and $\alpha\in(0,\frac{5}{4})$ (since $p_1=6\alpha-6$, $p_0=2\alpha-1$ and $\alpha_0=0$, see Lemma \ref{lema2.1} 2) and 3)).
Thereby,  by multiplying the inequality (\ref{wil5}) by $|\xi|^{p}$, and applying  Lemma \ref{lemaexponencial}, Lemma \ref{lema2.1} iii) and (\ref{P1}), it follows that
\begin{align*}
\nonumber|\xi|^{p}|\widehat{u}(\xi,t)|&\leq |\xi|^{p}|\widehat{u}_0(\xi)| + \int_{0}^t |\xi|^{p+1}e^{-  (t-\tau)|\xi|^{2\alpha}} \|u(\tau)\|_2^2d\tau\\
&\leq \|u_0\|_{\mathcal{Y}^p} + C\int_{0}^t (t-\tau)^{-\frac{p+1}{2\alpha}}(1+\tau)^{-\frac{3-2p_0}{2\alpha}}d\tau\\
&\leq \|u_0\|_{\mathcal{Y}^p} + C\int_{0}^t (t-\tau)^{-\frac{p+1}{2\alpha}}(1+\tau)^{-\frac{2\alpha-1-p}{2\alpha}}d\tau\\
&\leq \|u_0\|_{\mathcal{Y}^p} + C\int_{0}^t (t-\tau)^{-\frac{p+1}{2\alpha}}\tau^{-\frac{2\alpha-1-p}{2\alpha}}d\tau\\
&\leq \|u_0\|_{\mathcal{Y}^p} + C\,\beta(\tfrac{p+1}{2\alpha};\tfrac{2\alpha-1-p}{2\alpha}),
\end{align*}
for all $\xi\in \mathbb{R}^3$ and $t\geq0$ (because $p+1>0$, $2\alpha-1-p>0$, $\alpha>0$, $p\geq 6\alpha-6$ and $\tfrac{p+1}{2\alpha}+\tfrac{2\alpha-1-p}{2\alpha}=1$). Consequently, one infers
\begin{align*}
\|u\|_{L^\infty([0,\infty);\mathcal{Y}^p)}\leq \|u_0\|_{\mathcal{Y}^p} + C \,\beta(\tfrac{p+1}{2\alpha};\tfrac{2\alpha-1-p}{2\alpha}).
\end{align*}
(Recall   that $u_0\in \mathcal{Y}^p(\mathbb{R}^3)$). This gives ii) with $n=0.$\\\\
    iii) Assume that $0< \alpha< \frac{5}{6}$ (since $\alpha_0=0$ and $\alpha_1=\frac{5}{6}$, see Lemma \ref{lema2.1} 3)) and $p=-1$. Hence, by multiplying the inequality (\ref{wil5}) by $|\xi|^{-1}$ and applying   (\ref{P1}) and Lemma \ref{lema2.1} iii), we can write the following results:
\begin{align*}
\nonumber|\xi|^{-1}|\widehat{u}(\xi,t)|&\leq |\xi|^{-1}|\widehat{u}_0(\xi)| + \int_{0}^t e^{-  (t-\tau)|\xi|^{2\alpha}} \|u(\tau)\|_2^2d\tau\\
&\leq \|u_0\|_{\mathcal{Y}^{-1}} + C\int_{0}^t (1+\tau)^{-\frac{3-2p_0}{2\alpha}}d\tau\\
&\leq \|u_0\|_{\mathcal{Y}^{-1}} + \frac{2\alpha C}{3-2\alpha-2p_0}\\
&= \|u_0\|_{\mathcal{Y}^{-1}} + \frac{2\alpha C}{-(p_1+1)},
\end{align*}
for all $\xi\in \mathbb{R}^3$ and $t\geq0$ (recall that $-(p_1+1)=5-6\alpha>0$ (see Lemma \ref{lema2.1} 2)) and $0<\alpha<\frac{5}{6}$). Consequently, one infers
\begin{align*}
\|u\|_{L^\infty([0,\infty);\mathcal{Y}^{-1})}\leq \|u_0\|_{\mathcal{Y}^{-1}} + \frac{2\alpha C}{-(p_1+1)}.
\end{align*}
(Recall   that, in this case, $u_0\in \mathcal{Y}^{-1}(\mathbb{R}^3)$). This shows iii) with $n=0.$\\\\
\underline{2º Case}: Suppose that Proposition \ref{inducao} i)--iii) hold for $n-1\in \mathbb{N}\cup \{0\}$. Let us prove this same result, in the case $n\in \mathbb{N}$, as follows.\\\\
i) Consider that $p\leq p_{n}$ and $\alpha\in(\alpha_{n},\frac{5}{4})$. By using Lemma \ref{lema2.1}, one has
 \begin{align*}
 p\leq p_{n-1} \,\,  \hbox{  and  }\,\, \alpha_{n-1}< \alpha <\frac{5}{4},
 \end{align*}
 since $(p_n)_{n\geq0}$  and $(\alpha_n)_{n\geq0}$ are decreasing and increasing sequences, respectively. Thus, from the inductive  hypothesis, one concludes
    \begin{align}\label{Pn-1}
    \|u(t)\|_2^2\leq C(1+t)^{-\frac{3-2p_{n-1}}{2\alpha}},\quad \forall t>0.
    \end{align}
Thereby, by Lemma \ref{lemaexponencial} i) and (\ref{Pn-1}),  it follows that
\begin{align*}
\nonumber|\xi|^{p_n}|\widehat{u}(\xi,t)|&\leq |\xi|^{p_n-p}|\xi|^{p}|\widehat{u}_0(\xi)| + \int_{0}^t |\xi|^{p_n+1}e^{-  (t-\tau)|\xi|^{2\alpha}} \|u(\tau)\|_2^2d\tau\\
&\leq \|u_0\|_{\mathcal{Y}^p} + C\int_{0}^t (t-\tau)^{-\frac{p_n+1}{2\alpha}}(1+\tau)^{-\frac{3-2p_{n-1}}{2\alpha}}d\tau\\
&\leq \|u_0\|_{\mathcal{Y}^p} + C\,\beta(\tfrac{p_n+1}{2\alpha};\tfrac{3-2p_{n-1}}{2\alpha}),
\end{align*}
for all $\xi\in \mathbb{R}^3$ such that $|\xi|\leq1$ and $t\geq0$ (because $p\leq p_n$, $p_n+1>0$ (see Lemma \ref{lema2.1} vi)), $3-2p_{n-1}>0$ (see Lemma \ref{lema2.1} iv)), $\alpha>0$ (see Lemma \ref{lema2.1} v)) and $\tfrac{p_n+1}{2\alpha}+\tfrac{3-2p_{n-1}}{2\alpha}=1$ (see Lemma \ref{lema2.1} iii)). Therefore, we must have
\begin{align*}
\nonumber\sup_{t\in[0,\infty)}\{\sup_{|\xi|\leq1}\{|\xi|^{p_n}|\widehat{u}(\xi,t)|\}\}&\leq \|u_0\|_{\mathcal{Y}^p} + C\,\beta(\tfrac{p_n+1}{2\alpha};\tfrac{3-2p_{n-1}}{2\alpha}).
\end{align*}
(Recall   that $u_0\in \mathcal{Y}^p(\mathbb{R}^3)$). Moreover, Lemma \ref{lema2.1} iv) implies that $p_n<\frac{3}{2}$ and, as a result, Lemma \ref{lemaderivada}, Lemma \ref{lemaexponencial} ii) and (\ref{*12}) infer that
    \begin{align}\label{Pn}
    \|u(t)\|_2^2\leq C(1+t)^{-\frac{3-2p_{n}}{2\alpha}},\quad \forall t>0.
    \end{align}
This establishes i).\\\\
ii) Consider that $p_{n+1}\leq p< p_n$, $p>-1$ and $\alpha\in(\alpha_n,\frac{5}{4})$ to obtain, by (\ref{Pn}), Lemma \ref{lemaexponencial} i) and Lemma \ref{lema2.1}, that
\begin{align*}
\nonumber|\xi|^{p}|\widehat{u}(\xi,t)|&\leq |\xi|^{p}|\widehat{u}_0(\xi)| + \int_{0}^t |\xi|^{p+1}e^{-  (t-\tau)|\xi|^{2\alpha}} \|u(\tau)\|_2^2d\tau\\
&\leq \|u_0\|_{\mathcal{Y}^p} + C\int_{0}^t (t-\tau)^{-\frac{p+1}{2\alpha}}(1+\tau)^{-\frac{3-2p_n}{2\alpha}}d\tau\\
&\leq \|u_0\|_{\mathcal{Y}^p} + C\int_{0}^t (t-\tau)^{-\frac{p+1}{2\alpha}}(1+\tau)^{-\frac{2\alpha-1-p}{2\alpha}}d\tau\\
&\leq \|u_0\|_{\mathcal{Y}^p} + C\,\beta(\tfrac{p+1}{2\alpha};\tfrac{2\alpha-1-p}{2\alpha}),
\end{align*}
for all $\xi\in \mathbb{R}^3$ and $t\geq0$ (since $p+1>0$, $2\alpha-1> p_n>p$ (see Lemma \ref{lema2.1}), $\alpha>\alpha_n>0$  (see Lemma \ref{lema2.1} v)), $p\geq p_{n+1}$ and $\tfrac{p+1}{2\alpha}+\tfrac{2\alpha-1-p}{2\alpha}=1$). Hence, it follows that
\begin{align*}
\|u\|_{L^\infty([0,\infty);\mathcal{Y}^p)}\leq \|u_0\|_{\mathcal{Y}^p} + C\,\beta(\tfrac{p+1}{2\alpha};\tfrac{2\alpha-1-p}{2\alpha}).
\end{align*}
(Recall   that $u_0\in \mathcal{Y}^p(\mathbb{R}^3)$). This proves ii).\\\\
    iii) Assume that $\alpha_{n}< \alpha< \alpha_{n+1}$ and $p=-1$. Then, by Lemma \ref{lema2.1} vi), we conclude that
    \begin{align*}
    p=-1<p_n\,\, \hbox{ and }\,\, \alpha_n<\alpha<\frac{5}{4},
    \end{align*}
         since $\alpha\in(0,\frac{5}{4})$. Thus, by applying   (\ref{Pn}) and Lemma \ref{lema2.1} iii), one has
\begin{align*}
\nonumber|\xi|^{-1}|\widehat{u}(\xi,t)|&\leq |\xi|^{-1}|\widehat{u}_0(\xi)| + \int_{0}^t e^{-  (t-\tau)|\xi|^{2\alpha}} \|u(\tau)\|_2^2d\tau\\
&\leq \|u_0\|_{\mathcal{Y}^{-1}} + C\int_{0}^t (1+\tau)^{-\frac{3-2p_n}{2\alpha}}d\tau\\
&\leq \|u_0\|_{\mathcal{Y}^{-1}} + \frac{2\alpha C}{3-2\alpha-2p_n}\\
&= \|u_0\|_{\mathcal{Y}^{-1}} + \frac{2\alpha C}{-(p_{n+1}+1)},
\end{align*}
for all $\xi\in \mathbb{R}^3$ and $t\geq0$ (recall that $p_{n+1}+1<0$ (see Lemma \ref{lema2.1} vii)) and $\alpha>\alpha_n>0$ (see Lemma \ref{lema2.1} v)). As a consequence, one infers
\begin{align*}
\|u\|_{L^\infty([0,\infty);\mathcal{Y}^{-1})}\leq \|u_0\|_{\mathcal{Y}^{-1}} + \frac{2\alpha C}{-(p_{n+1}+1)}.
\end{align*}
(Recall   that, in this case, $u_0\in \mathcal{Y}^{-1}(\mathbb{R}^3)$). This shows iii).\\\\

These arguments above prove that Proposition \ref{inducao} i), ii) and iii) hold for any $n\in\mathbb{N}\cup\{0\}.$

Lastly, consider that $p\geq2\alpha-1$.  By multiplying the inequality (\ref{wil5}) by $|\xi|^{p}$, and applying (\ref{*121}), it follows that
\begin{align}\label{wn1}
\nonumber|\xi|^{p}|\widehat{u}(\xi,t)|&\leq |\xi|^{p}|\widehat{u}_0(\xi)| + \int_{0}^t |\xi|^{p+1}e^{-  (t-\tau)|\xi|^{2\alpha}} \|u(\tau)\|_2^2d\tau\\
&\leq \|u_0\|_{\mathcal{Y}^{p}} + \|u_0\|_2^2\int_{0}^t |\xi|^{p-(2\alpha-1)}|\xi|^{2\alpha}e^{-  (t-\tau)|\xi|^{2\alpha}}d\tau\\
\nonumber&\leq \|u_0\|_{\mathcal{Y}^{p}} + \|u_0\|_2^2\int_{0}^t |\xi|^{2\alpha}e^{-  (t-\tau)|\xi|^{2\alpha}}d\tau\\
\nonumber&\leq \|u_0\|_{\mathcal{Y}^{p}} + \|u_0\|_2^2,
\end{align}
for all $\xi\in \mathbb{R}^3$ such that $|\xi|\leq 1$ and $t\geq0$ (since $p-(2\alpha-1)\geq0)$. Therefore, it is always true that
\begin{align*}
\sup_{t\in[0,\infty)}\{\sup_{|\xi|\leq 1} \{|\xi|^p|\hat{u}(\xi,t)|\}\}\leq \|u_0\|_{\mathcal{Y}^{p}} + \|u_0\|_2^2.
\end{align*}
(Recall   that  $u_0\in L^{2}(\mathbb{R}^3)\cap\mathcal{Y}^{p}(\mathbb{R}^3)$). This establishes the proof of (\ref{2alpha-1}).

Also, if $p=2\alpha-1$, by (\ref{wn1}), we can write
\begin{align*}
\|u\|_{L^\infty([0,\infty);\mathcal{Y}^{2\alpha-1})}\leq \|u_0\|_{\mathcal{Y}^{2\alpha-1}} + \|u_0\|_2^2.
\end{align*}
(Recall   that, in this case,  $u_0\in L^{2}(\mathbb{R}^3)\cap\mathcal{Y}^{2\alpha-1}(\mathbb{R}^3)$). This proves  (\ref{wn}).

\end{proof}

\hspace{-0.6cm}\underline{\textbf{Proof of Theorem \ref{teoremaexistenciaw1}}}:\\

Now, we are ready to establish a proof of the most important result of this work: Theorem \ref{teoremaexistenciaw1}. Let us point out that   our arguments presented below were motivated by the papers \cite{novo,nonlinear}.

We shall split our prove into two cases.\\\\
\underline{1º Case}: Assume that $\alpha\in \cup_{i=0}^{m_0+1} (\alpha_i,\alpha_{i+1})$, where $m_0$ has been found in (\ref{sm0}) (see also (\ref{alphaepm0})).\\\\
At first, it is worth to recall that (\ref{alphaepm0}) informs that $0<\alpha< \alpha_{m_0+2}$. Secondly, in this case, it follows that
\begin{align}\label{wil6}
\alpha\in (\alpha_{i_0},\alpha_{i_0+1}),\quad \hbox{for some } i_0=0,1,2,...,m_0+1.
\end{align}
Consequently, by Lemma \ref{lema2.1} v), vi) and vii), one reaches
\begin{align}\label{pi0alphai0}
p_{i_0 +1}<-1<p_{i_0} \,\, \hbox{ and } \,\, \alpha_{i_0+1}<\frac{5}{4}.
\end{align}
$\bullet$ Consider that $p=-1.$\\\\
From (\ref{wil6}) and Proposition \ref{inducao} iii), we can write
\begin{align}\label{wil8}
\|u\|_{L^\infty([0,\infty);\mathcal{Y}^{-1})}\leq \|u_0\|_{\mathcal{Y}^{-1}} + \frac{2\alpha C}{-(p_{i_0+1}+1)}.
\end{align}
It is important to point out that $p_{i_0 +1}+1<0$ (see (\ref{pi0alphai0})).\\\\
$\bullet$ Consider that $p\in (-1,p_{i_0}).$\\\\
In this case, by applying (\ref{pi0alphai0}) and (\ref{wil6}), we can conclude that
\begin{align*}
p_{i_0 +1}<p<p_{i_0}, \,\,p>-1 \,\, \hbox{ and } \,\, \alpha_{i_0}<\alpha<\frac{5}{4}.
\end{align*}
Therefore, Proposition \ref{inducao} ii) implies that
\begin{align}\label{wil10}
\|u\|_{L^\infty([0,\infty);\mathcal{Y}^p)}\leq \|u_0\|_{\mathcal{Y}^p} + C\,\beta(\tfrac{p+1}{2\alpha};\tfrac{2\alpha-1-p}{2\alpha}).
\end{align}
$\bullet$ Consider that $p\in [p_{i_0},p_0)=\cup_{k=1}^{i_0}[p_k,p_{k-1}).$\\\\
First of all, observe that if $i_0=0$; then, the two previous cases would prove that
\begin{align*}
\|u\|_{L^\infty([0,\infty);\mathcal{Y}^p)}\leq C,\quad\forall p\in[-1,2\alpha-1),
\end{align*}
since $u_0\in \mathcal{Y}^p(\mathbb{R}^3)$. Thus, we shall assume  $i_0\geq1$. Thereby, it is true that
\begin{align}\label{wil7}
p_{k_0}\leq p <p_{k_0 -1},\quad \hbox{for some } k_0=1,2,...,i_0.
\end{align}
Moreover, by Lemma \ref{lema2.1}, (\ref{wil6}), (\ref{pi0alphai0}) and (\ref{wil7}), one reaches
\begin{align}\label{pk0alphak0}
-1<p \,\, \hbox{ and } \,\, \alpha_{k_0-1}<\alpha<\frac{5}{4},
\end{align}
since $1\leq k_0\leq i_0$. As a result, (\ref{wil7}), (\ref{pk0alphak0}) and Proposition \ref{inducao} ii) imply that
\begin{align}\label{wil11}
\|u\|_{L^\infty([0,\infty);\mathcal{Y}^p)}\leq \|u_0\|_{\mathcal{Y}^p} + C\,\beta(\tfrac{p+1}{2\alpha};\tfrac{2\alpha-1-p}{2\alpha}).
\end{align}
$\bullet$ Consider that $p\in[p_0,\frac{3}{2})$ (recall that $p_0=2\alpha-1$).\\\\
By using (\ref{2alpha-1}), we obtain
\begin{align}\label{wil9}
\sup_{t\in[0,\infty)}\{\sup_{|\xi|\leq 1} \{|\xi|^p|\hat{u}(\xi,t)|\}\}\leq \|u_0\|_{\mathcal{Y}^{p}} + \|u_0\|_2^2.
\end{align}
Moreover, for $p=p_0 (=2\alpha-1)$, by (\ref{wn}), we deduce that
\begin{align}\label{wn2}
\|u\|_{L^\infty([0,\infty);\mathcal{Y}^{2\alpha-1})}\leq \|u_0\|_{\mathcal{Y}^{2\alpha-1}} + \|u_0\|_2^2.
\end{align}

Lastly, (\ref{wil8}), (\ref{wil10}), (\ref{wil11}) and (\ref{wil9}) show that
\begin{align}\label{wil12}
\sup_{t\in[0,\infty)}\{\sup_{|\xi|\leq 1} \{|\xi|^p|\hat{u}(\xi,t)|\}\}\leq C,\quad\forall p\in[-1,\tfrac{3}{2}),
\end{align}
since $u_0\in L^2(\mathbb{R}^3)\cap \mathcal{Y}^p(\mathbb{R}^3)$. It is also important to emphasize that $$u\in L^\infty([0,\infty);\mathcal{Y}^p(\mathbb{R}^3)),\quad\forall p\in[-1,2\alpha-1],$$ by (\ref{wil8}), (\ref{wil10}), (\ref{wil11}) and (\ref{wn2}). This proves (\ref{wn10}).

Therefore, by  (\ref{*12}), (\ref{wil12}) (notice that $p<\frac{3}{2}$) and Lemma \ref{lemaderivada}, one deduces
  \begin{align*}
    \|u(t)\|_2^2\leq C(1+t)^{-\frac{3-2p}{2\alpha}},\quad \forall t>0,
    \end{align*}
    where $p\in[-1,\tfrac{3}{2}).$ This establishes (\ref{l2}).\\\\
\underline{2º Case}: Assume that $\alpha= \alpha_i$ with $i=1,2,...,m_0+1$,  where $m_0$ has been found in (\ref{sm0}) (see also (\ref{alphaepm0})).\\\\
By observing Lemma \ref{lema2.1} v) and vi), we have
\begin{align}\label{wil22}
\alpha_{i-1}< \alpha<\frac{5}{4}\,\, \hbox{ and }\,\, p_0\geq p_{i-1}>-1.
\end{align}
$\bullet$ Consider that $p=-1$.\\\\
By adding that $u_0\in L^1(\mathbb{R}^3)$ (this is necessary only in this case of our proof) as well, and using (\ref{nonlinear1}) (see \cite{nonlinear}), we obtain
\begin{align*}
\nonumber|\xi|^{-1}|\widehat{u}(\xi,t)|&\leq |\xi|^{-1}|\widehat{u}_0(\xi)| + \int_{0}^t e^{-  (t-\tau)|\xi|^{2\alpha}} \|u(\tau)\|_2^2d\tau\\
&\leq \|u_0\|_{\mathcal{Y}^{-1}} + C\int_{0}^t (1+\tau)^{-\frac{3}{2\alpha}}d\tau\\
&\leq \|u_0\|_{\mathcal{Y}^{-1}} + \frac{2\alpha C}{3-2\alpha},
\end{align*}
for all $\xi\in \mathbb{R}^3$ and $t\geq0$ (recall that $3-2\alpha>0$  (see (\ref{wil22})) and $\alpha>0$ (see Lemma \ref{lema2.1} v)). As a consequence, one infers
\begin{align}\label{wil30}
\|u\|_{L^\infty([0,\infty);\mathcal{Y}^{-1})}\leq \|u_0\|_{\mathcal{Y}^{-1}} + \frac{2\alpha C}{3-2\alpha}.
\end{align}
$\bullet$ Consider that $p\in (-1,p_{i-1}).$\\\\
Lemma \ref{lema2.1} viii) implies that $p_i=-1$ (since $\alpha=\alpha_i$). As a consequence of (\ref{wil22}), we deduce that
\begin{align}\label{wil20}
p_i<p<p_{i-1},\,p>-1\,\,  \hbox{ and }\,\, \alpha_{i-1}<\alpha<\frac{5}{4}.
\end{align}
Thereby, from (\ref{wil20}) and Proposition \ref{inducao} ii), it follows that
\begin{align}\label{wil32}
\|u\|_{L^\infty([0,\infty);\mathcal{Y}^p)}\leq \|u_0\|_{\mathcal{Y}^p} + C\,\beta(\tfrac{p+1}{2\alpha};\tfrac{2\alpha-1-p}{2\alpha}).
\end{align}
$\bullet$ Consider that $p\in [p_{i-1},p_0)=\cup_{l=1}^{i-1}[p_l,p_{l-1}).$\\\\
At first, observe that if $i=1$; then, the two previous cases would prove that
\begin{align*}
\|u\|_{L^\infty([0,\infty);\mathcal{Y}^p)}\leq C,\quad\forall p\in[-1,2\alpha-1),
\end{align*}
since $u_0\in \mathcal{Y}^p(\mathbb{R}^3)$. Hence, we shall assume  $i\geq2$. Consequently, in this case, it holds that
\begin{align}\label{wil21}
p_{l_0}\leq p <p_{l_0 -1},\quad \hbox{for some } l_0=1,2,...,i-1.
\end{align}
In addition, by Lemma \ref{lema2.1}, (\ref{wil22}) and (\ref{wil21}), one reaches
\begin{align}\label{pl0alphal0}
p>-1 \,\, \hbox{ and } \,\, \alpha_{l_0-1}<\alpha<\frac{5}{4},
\end{align}
since $1\leq l_0\leq i-1$. As a result, (\ref{wil21}), (\ref{pl0alphal0}) and Proposition \ref{inducao} ii) imply that
\begin{align}\label{wil33}
\|u\|_{L^\infty([0,\infty);\mathcal{Y}^p)}\leq \|u_0\|_{\mathcal{Y}^p} + C\,\beta(\tfrac{p+1}{2\alpha};\tfrac{2\alpha-1-p}{2\alpha}).
\end{align}
$\bullet$ Consider that $p\in[p_0,\frac{3}{2})$ (recall that $p_0=2\alpha-1$).\\\\
By using (\ref{2alpha-1}), we obtain
\begin{align}\label{wil31}
\sup_{t\in[0,\infty)}\{\sup_{|\xi|\leq 1} \{|\xi|^p|\hat{u}(\xi,t)|\}\}\leq \|u_0\|_{\mathcal{Y}^{p}} + \|u_0\|_2^2.
\end{align}
Moreover, for $p=p_0 (=2\alpha-1)$, by (\ref{wn}), we deduce that
\begin{align}\label{wn5}
\|u\|_{L^\infty([0,\infty);\mathcal{Y}^{2\alpha-1})}\leq \|u_0\|_{\mathcal{Y}^{2\alpha-1}} + \|u_0\|_2^2.
\end{align}

At last, (\ref{wil30}),  (\ref{wil32}), (\ref{wil33}) and (\ref{wil31}) show that
\begin{align}\label{wil34}
\sup_{t\in[0,\infty)}\{\sup_{|\xi|\leq 1} \{|\xi|^p|\hat{u}(\xi,t)|\}\}\leq C,\quad\forall p\in[-1,\tfrac{3}{2}),
\end{align}
since $u_0\in L^2(\mathbb{R}^3)\cap \mathcal{Y}^p(\mathbb{R}^3)$. These arguments also imply that
$$u \in L^\infty([0,\infty);\mathcal{Y}^p(\mathbb{R}^3)), \quad\forall p\in[-1,2\alpha-1],$$
by (\ref{wil30}),  (\ref{wil32}), (\ref{wil33}) and (\ref{wn5}). This proves (\ref{wn10}).

Thus, by using  (\ref{*12}), (\ref{wil34}) (notice that $p<\frac{3}{2}$) and Lemma \ref{lemaderivada}, we reach
  \begin{align*}
    \|u(t)\|_2^2\leq C(1+t)^{-\frac{3-2p}{2\alpha}},\quad \forall t>0,
    \end{align*}
where $p\in[-1,\tfrac{3}{2})$. This establishes (\ref{l2}). Therefore,  Theorem \ref{teoremaexistenciaw1} is proved.

\caixa

\hspace{-0.6cm}\underline{\textbf{Proof of Corollary \ref{corolario}}}:\\

In order to finish this work, we shall establish the proof of Corollary \ref{corolario}. We shall adapt the arguments presented by \cite{novo}.

First of all, fix $t\geq0$ and notice that
\begin{align}\label{wc1}
\nonumber\|u(t)\|_{\dot{H}^{-\delta}}^2&=\int_{|\xi|\leq N(t)} |\xi|^{-2\delta} |\widehat{u}(t)|^2 d\xi+\int_{|\xi|> N(t)} |\xi|^{-2\delta} |\widehat{u}(t)|^2 d\xi\\
\nonumber&=\int_{|\xi|\leq N(t)} |\xi|^{-2\delta} |\widehat{u}(t)|^2 d\xi+ [N(t)]^{-2\delta}\int_{|\xi|> N(t)}  |\widehat{u}(t)|^2 d\xi\\
&\leq\int_{|\xi|\leq N(t)} |\xi|^{-2\delta} |\widehat{u}(t)|^2 d\xi+ C[N(t)]^{-2\delta}\|u(t)\|_2^2,
\end{align}
by using Plancherel's identity and the fact that $\delta\geq0$.
\\\\
i) Assume that $p\in[-1,2\alpha-1]$.\\\\
Let $N(t)$ be the following  real number:
\begin{align}\label{wn6}
N(t)=\|u(t)\|_2^{\frac{2}{3-2p}}\|u(t)\|_{\mathcal{Y}^p}^{-\frac{2}{3-2p}}.
\end{align}
On the other hand,  it is true that
\begin{align}\label{wc2}
\nonumber\int_{|\xi|\leq N(t)} |\xi|^{-2\delta} |\widehat{u}(t)|^2 d\xi&= \int_{|\xi|\leq N(t)} |\xi|^{-2\delta-2p} |\xi|^{2p}|\widehat{u}(t)|^2 d\xi\\
\nonumber&\leq \|u(t)\|_{\mathcal{Y}^p}^2 \int_{|\xi|\leq N(t)} |\xi|^{-2\delta-2p}  d\xi\\
\nonumber&= C\|u(t)\|_{\mathcal{Y}^p}^2 \int_{0}^{N(t)} r^{2-2\delta-2p}  dr\\
&= C\|u(t)\|_{\mathcal{Y}^p}^2 N(t)^{3-2\delta-2p},
\end{align}
since $3-2\delta-2p>0$. By replacing (\ref{wc2}) in (\ref{wc1}), one infers
\begin{align}\label{wc3}
\|u(t)\|_{\dot{H}^{-\delta}}^2&\leq C[\|u(t)\|_{\mathcal{Y}^p}^2 N(t)^{3-2\delta-2p}+ N(t)^{-2\delta}\|u(t)\|_2^2],
\end{align}
for all $t\geq0$. By observing (\ref{wc3}), (\ref{wn6}), (\ref{wil8}), (\ref{wil10}), (\ref{wil11}), (\ref{wn2}), (\ref{wil30}),  (\ref{wil32}), (\ref{wil33}) and (\ref{wn5}) (recall that $p\in[-1,2\alpha-1]$), we can write
\begin{align*}
\nonumber\|u(t)\|_{\dot{H}^{-\delta}}^2&\leq C\|u(t)\|_{\mathcal{Y}^p}^{\frac{4\delta}{3-2p}} \|u(t)\|_2^{\frac{2(3-2\delta-2p)}{3-2p}}\leq C \|u(t)\|_2^{\frac{2(3-2\delta-2p)}{3-2p}},
\end{align*}
for all $t\geq0$, since $3-2\delta-2p>0$ and $\delta\geq0$. By applying (\ref{l2}), one reaches
\begin{align*}
\nonumber\|u(t)\|_{\dot{H}^{-\delta}}^2\leq C (1+t)^{-\frac{3-2\delta-2p}{2\alpha}}, \quad \forall t>0,
\end{align*}
whether $p\in[-1,2\alpha-1]$. This proves (\ref{hdelta1}).
\\\\
ii) Assume that $p\in[2\alpha-1,\frac{3}{2})$.\\\\
Let $N(t)$ be the following  real number:
\begin{align}\label{wfn6}
N(t)=\|u(t)\|_2^{\frac{2}{3-2p}}.
\end{align}
On the other hand, by (\ref{wn1}), we conclude
\begin{align}\label{wfn1}
\nonumber|\xi|^{p}|\widehat{u}(\xi,t)|&\leq \|u_0\|_{\mathcal{Y}^{p}} + \|u_0\|_2^2[N(t)]^{p-(2\alpha-1)}\int_{0}^t |\xi|^{2\alpha}e^{-  (t-\tau)|\xi|^{2\alpha}}d\tau\\
&\leq \|u_0\|_{\mathcal{Y}^{p}} + \|u_0\|_2^2[N(t)]^{p-(2\alpha-1)},
\end{align}
for all $\xi\in \mathbb{R}^3$ such that $|\xi|\leq N(t)$ (since $p-(2\alpha-1)\geq0)$. Furthermore, (\ref{wfn1}) implies that
\begin{align}\label{wfc2}
\nonumber\int_{|\xi|\leq N(t)} |\xi|^{-2\delta} |\widehat{u}(t)|^2 d\xi&= \int_{|\xi|\leq N(t)} |\xi|^{-2\delta-2p} |\xi|^{2p}|\widehat{u}(t)|^2 d\xi\\
\nonumber&\leq C[1+N(t)^{p-2\alpha+1}]^2 \int_{|\xi|\leq N(t)} |\xi|^{-2\delta-2p}  d\xi\\
\nonumber&\leq C[1+N(t)^{2p-4\alpha+2}] \int_{0}^{N(t)} r^{2-2\delta-2p}  dr\\
&= C[N(t)^{3-2\delta-2p}+N(t)^{5-4\alpha-2\delta}] ,
\end{align}
since $3-2\delta-2p>0$ and $u_0\in L^2(\mathbb{R}^3)\cap \mathcal{Y}^p(\mathbb{R}^3) $. By replacing (\ref{wfc2}) in (\ref{wc1}), one infers
\begin{align}\label{wfc3}
\|u(t)\|_{\dot{H}^{-\delta}}^2&\leq C[N(t)^{3-2\delta-2p}+N(t)^{5-4\alpha-2\delta}+ N(t)^{-2\delta}\|u(t)\|_2^2],
\end{align}
for all $t\geq0$. By  (\ref{wfc3}) and (\ref{wfn6}), it follows that
\begin{align*}
\nonumber\|u(t)\|_{\dot{H}^{-\delta}}^2&\leq C[\|u(t)\|_2^{\frac{2(3-2\delta-2p)}{3-2p}}+\|u(t)\|_2^{\frac{2(5-4\alpha-2\delta)}{3-2p}}],
\end{align*}
for all $t\geq0$. Since $3-2\delta-2p>0$, $p\in[2\alpha-1,\frac{3}{2})$ and $\delta\geq0$, by applying (\ref{l2}), we deduce
\begin{align*}
\nonumber\|u(t)\|_{\dot{H}^{-\delta}}^2&\leq C[(1+t)^{-\frac{3-2\delta-2p}{2\alpha}}+(1+t)^{-\frac{5-4\alpha-2\delta}{2\alpha}}],
\end{align*}
 for all $t>0$. 
This establishes (\ref{hdelta}).

These two items prove Corollary \ref{corolario}.

\caixa

\subsection*{Acknowledgments}
Wilberclay G. Melo is partially supported by CNPq grant 309880/2021-1.

\end{document}